\documentclass[a4paper]{amsart}

\usepackage{amssymb}
\usepackage[arrow]{xy}
\usepackage{pb-diagram,pb-xy}
\usepackage{graphicx,color,float}

\theoremstyle{plain}
\newtheorem{theorem}{Theorem}[section]
\newtheorem{proposition}[theorem]{Proposition}
\newtheorem{corollary}[theorem]{Corollary}
\newtheorem{lemma}[theorem]{Lemma}

\theoremstyle{definition}
\newtheorem{definition}[theorem]{Definition}

\def\Z{\mathbb{Z}}
\def\Q{\mathbb{Q}}

\def\C{\mathbb{C}}

\def\K{\mathcal{K}}

\def\Ker{\operatorname{Ker}}
\def\Im{\operatorname{Im}}
\def\Hom{\operatorname{Hom}}

\def\rank{\operatorname{rank}}

\def\G1{\Gamma^{(1)}_r}
\def\bbk{\mathbb{K}}

\def\zp{\mathbb{Z}_{(p)}}
\def\z2{\mathbb{Z}_{(2)}}
\def\A{\mathcal{A}}
\def\F{\mathcal{F}}
\def\C{\mathcal{C}}

\def\to{\mathchoice{\longrightarrow}{\rightarrow}{\rightarrow}{\rightarrow}}
\makeatletter
\newcommand{\shortxra}[2][]{\ext@arrow 0359\rightarrowfill@{#1}{#2}}
\def\longrightarrowfill@{\arrowfill@\relbar\relbar\longrightarrow}
\newcommand{\longxra}[2][]{\ext@arrow 0359\longrightarrowfill@{#1}{#2}}
\renewcommand{\xrightarrow}[2][]{\mathchoice{\longxra[#1]{#2}}%
  {\shortxra[#1]{#2}}{\shortxra[#1]{#2}}{\shortxra[#1]{#2}}}
\makeatother

\makeatletter
\def\Nopagebreak{\@nobreaktrue\nopagebreak}
\makeatother

\let\ifshowaddedremoved=\iftrue

\newtoks\trashtoks
\long\def\totrash#1\endremoved{}
\ifshowaddedremoved
  \def\beginremoved{\begin{color}{blue}\textsf{[Removed]}\space}
  \def\endremoved{\end{color}}
  \def\beginadded{\begin{color}{red}\textsf{[Added]}\space}
  \def\endadded{\end{color}}
\else
  \def\beginremoved{\totrash}
  \def\endremoved{}
  \def\beginadded{}
  \def\endadded{}
\fi

\begin{document}

\title
{Covering link calculus and iterated Bing doubles}

\author{Jae Choon Cha}

\address{Department of Mathematics and Pohang Mathematics Institute \\
  Pohang University of Science and Technology\\ Pohang, Gyungbuk
  790--784 \\ Republic of Korea}

\email{jccha@postech.ac.kr}

\author{Taehee Kim}

\address{Department of Mathematics \\ Konkuk University \\ Seoul 143--701 \\
  Republic of Korea}

\email{tkim@konkuk.ac.kr}


\def\subjclassname{\textup{2000} Mathematics Subject Classification}
\expandafter\let\csname subjclassname@1991\endcsname=\subjclassname
\expandafter\let\csname subjclassname@2000\endcsname=\subjclassname


%


\begin{abstract}
  We give a new geometric obstruction to the iterated Bing double of a
  knot being a slice link: for $n>1$ the $(n+1)$st iterated Bing
  double of a knot is rationally slice if and only if the $n$th
  iterated Bing double of the knot is rationally slice.  The main
  technique of the proof is a covering link construction simplifying a
  given link.  We prove certain similar geometric obstructions for
  $n\le 1$ as well.  Our results are sharp enough to conclude, when
  combined with algebraic invariants, that if the $n$th iterated Bing
  double of a knot is slice for some $n$, then the knot is
  algebraically slice.  Also our geometric arguments applied to the
  smooth case show that the Ozsv\'ath--Szab\'o and Manolescu--Owens
  invariants give obstructions to iterated Bing doubles being slice.
  These results generalize recent results of Harvey, Teichner,
  Cimasoni, Cha and Cha--Livingston--Ruberman.  As another application,
  we give explicit examples of algebraically slice knots with
  non-slice iterated Bing doubles by considering von Neumann
  $\rho$--invariants and rational knot concordance.  Refined versions
  of such examples are given, that take into account the
  Cochran--Orr--Teichner filtration.
\end{abstract}

\maketitle

\section{Introduction}

\begin{figure}[ht]
  \begin{center}
    \includegraphics[scale=.9]{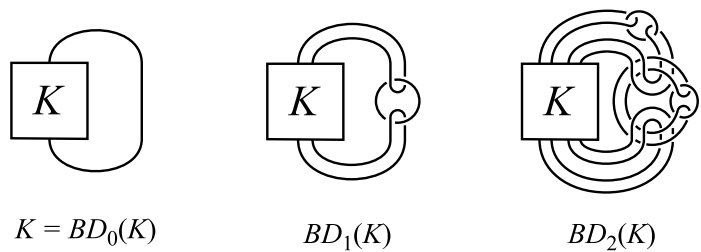}
  \end{center}
  \caption{Iterated Bing doubles $BD_n(K)$ of a knot
    $K$}\label{fig:iterated-bing-double-example.eps}
\end{figure}

The \emph{Bing double} $BD(K)$ of a knot $K$ is defined to be the
2--component link obtained by taking two zero-linking parallel copies
of $K$ and introducing positive and negative clasps, as in
Figure~\ref{fig:iterated-bing-double-example.eps}.  Taking the Bing
double of each component of $BD(K)$, we obtain the second iterated
Bing double $BD_2(K)$ of $K$.  Iterating this process, we define the
\emph{$n$th iterated Bing double} $BD_n(K)$ of $K$, which is a link
with $2^n$ components.  As our convention, for $n=0$, $BD_n(K)$
designates $K$ itself.

The problem of deciding whether $BD_n(K)$ is slice for some $n\ge 1$
has been studied actively, partly motivated by the relationship with
the 4--dimensional surgery theory.  We recall that a link $L$ in the
3--sphere $S^3$ is a \emph{slice link} if the components of $L$ bound
disjoint locally flat 2--disks in the 4--ball~$B^4$.  A fact that makes
the problem more interesting is that many previously known
obstructions to being a slice link vanish for any (iterated) Bing
double.  For an excellent discussion on this, the reader is referred
to Cimasoni's paper~\cite{Cimasoni:2006-1}.

It can be seen easily that if $BD_n(K)$ is slice, then so is
$BD_{n+1}(K)$.  Consequently if $K$ is slice then all iterated Bing
doubles of $K$ are slice.  The converse is a well known open problem.
Recently, there has been significant progress that enables us to
extract obstructions for (iterated) Bing doubles to be slice, and
consequently partial results on the converse.

A first remarkable result in this direction has been proved by
Harvey~\cite{Harvey:2006-1}, and Teichner (unpublished),
independently, using von Neumann $\rho$--invariants: if $BD_n(K)$ is
slice for some $n$, then the integral of the Levine--Tristram signature
of $K$ over the unit circle is zero.  In \cite{Cimasoni:2006-1}
Cimasoni proved that $K$ is algebraically slice if $BD(K)$ is
``boundary'' slice in the sense of
\cite{Cappell-Shaneson:1980-1,Ko:1987-1,Mio:1987-1,Duval:1986-1}.

As an application of his Hirzebruch-type intersection form defect
invariants, the first author found a new technique to detect non-slice
iterated Bing doubles which is effective even for knots of finite
order in the knot concordance group~\cite{Cha:2007-5}.  Using this, he
generalized the result of Harvey and Teichner by proving that for any
$n$ the Levine--Tristram signature function of $K$ is determined by
(the concordance class of) $BD_n(K)$, and also found infinitely many
amphichiral knots with non-slice iterated Bing doubles \cite[Theorems
1.5 and 1.6]{Cha:2007-5}.  In particular he gave the first proof that
any iterated Bing double of the figure eight knot is not slice.
Subsequent to this, Livingston, Ruberman, and the first author proved
that if $BD(K)$ is slice, then $K$ is algebraically
slice~\cite[Theorem 1]{Cha-Livingston-Ruberman:2006-1}.  Recently,
Cochran, Harvey, and Leidy showed that there are algebraically slice
knots with non-slice iterated Bing doubles using higher-order
$L^2$--signatures \cite{Cochran-Harvey-Leidy:2007-1}.

In this paper, we extend the aforementioned results on slicing
iterated Bing doubles.  First we prove a geometric result that the
converse of the fact ``$BD_n(K)$ is slice $\Longrightarrow$
$BD_{n+1}(K)$ is slice'' is \emph{rationally} true for higher~$n$:

\begin{theorem}\label{thm:main1}
  For any $n>1$, $BD_{n+1}(K)$ is rationally slice if and only if
  $BD_{n}(K)$ is rationally slice.
\end{theorem}

Here, as in \cite{Cha-Ko:2000-1,Cha:2003-1}, a link $L$ is said to be
a \emph{rationally slice link} if its ambient space is the boundary of
some rational homology 4--ball $W$ and there are disjoint locally flat
2--disks in~$W$ with boundary~$L$.  For a prime $p$, a
\emph{$\zp$--slice link} is defined similarly, namely slicing disks
exist in a $\zp$--homology ball instead.  (Here $\zp$ denotes the
localization of $\Z$ at the prime~$p$.)  A slice link is $\zp$--slice
for every prime $p$.  A link is $\zp$--slice for some prime $p$ if and
only if it is rationally slice.

In fact, we prove the $\zp$--analogue of Theorem~\ref{thm:main1}, from
which Theorem~\ref{thm:main1} follows immediately.  As the main
technique of the proof, we perform certain \emph{iterated covering
  link calculus} for iterated Bing doubles.  Given a link $L$ in a
$\zp$--homology sphere, the $p^a$--fold cyclic cover of the ambient
space branched over a component of $L$ becomes another $\zp$--homology
sphere and the pre-image of $L$ can be regarded as a new link.  The
idea of taking such a branched cover was first applied to
(non-iterated) Bing doubles in the work of
Cha--Livingston--Ruberman~\cite{Cha-Livingston-Ruberman:2006-1}.  We
perform a more sophisticated covering link calculus, by iterating the
process of taking branched coverings and taking sublinks; we call
links obtained in this way \emph{$p$--covering links}.  (See
Section~\ref{sec:covering}.)  The essential part of the proof of
Theorem~\ref{thm:main1} is the following: for $n>1$, $BD_{n}(K)$ is a
$p$--covering link of a more complicated link, namely $BD_{n+1}(K)$.
(See Proposition~\ref{prop:covering}.)

For the case of $n\le 1$, we do not know whether or not $BD_n(K)$ is a
$p$--covering link of $BD_{n+1}(K)$.  However, similarly to results for
$n=0$ in \cite{Cha-Livingston-Ruberman:2006-1}, our iterated covering
link technique can be used to show that certain band sums of (parallel
copies of) $K$ and its reverse $K^r$ are $\zp$--slice if $BD_{n+1}(K)$
is $\zp$--slice for $n\le1$.  (For example, see
Proposition~\ref{prop:connected sum} and its use in
Section~\ref{sec:algebraic}.)  The following result is a simple
special case:

\begin{proposition}\label{prop:reverse}
  If $BD_{n}(K)$ is $\zp$--slice for some $n\ge 0$, then $2K\#2K^r$ is
  $\zp$--slice.
\end{proposition}

We remark that our covering link calculus argument works in both
topological and smooth cases, so that Theorem~\ref{thm:main1} and
Proposition~\ref{prop:reverse} hold in the smooth case as well.

Combining our geometric results with previously known facts on
algebraic invariants of the $\Z_{(2)}$-concordance
group~\cite{Cochran-Orr:1993-1,Cha:2003-1}, we can deduce the
following second main theorem of this paper:

\begin{theorem}\label{thm:main2}
  For any $n$, if $BD_n(K)$ is slice, then $K$ is algebraically slice.
\end{theorem}

This generalizes the result for $BD_1(K)$ due to
Cha--Livingston--Ruberman
\cite[Theorem~1]{Cha-Livingston-Ruberman:2006-1} and generalizes the
first author's Levine--Tristram signature obstruction for $BD_n(K)$ to
be a slice link~\cite{Cha:2007-5}.  Theorem~\ref{thm:main2} can also
be used to show the following result which was first shown
in~\cite{Cha:2007-5}: there exist infinitely many knots $K$ such that
$K$ is amphichiral (so that it has order 2 in the knot concordance
group) but $BD_n(K)$ is not slice for any~$n$.

\begin{figure}[H]
  \begin{center}
    \includegraphics[scale=.95]{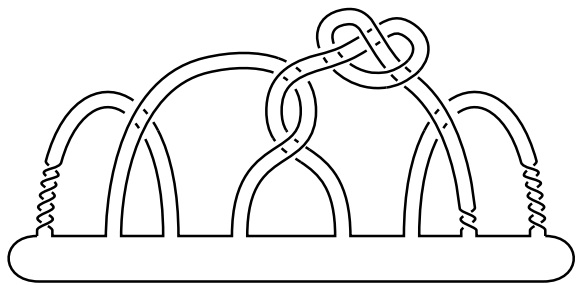}
  \end{center}
  \caption{An algebraically slice knot $K$ with $BD_n(K)$
    nonslice}\label{fig:alg-slice-with-nonslice-bd}
\end{figure}

Our geometric results can also be applied to investigate
(non-)sliceness of iterated Bing doubles of algebraically slice knots.
Recently Cochran--Harvey--Leidy \cite{Cochran-Harvey-Leidy:2007-1}
showed the existence of algebraically slice knots $K$ with non-slice
$BD_n(K)$.  In this paper, using techniques different from the ones in
\cite{Cochran-Harvey-Leidy:2007-1}, we construct \emph{explicit}
examples:

\begin{theorem}\label{thm:main3}
  The knot $K$ illustrated in
  Figure~\ref{fig:alg-slice-with-nonslice-bd} is algebraically slice
  but $BD_n(K)$ is not slice for any~$n$.
\end{theorem}

In fact, our method gives infinitely many explicit examples.  For
example, for any odd prime $q$, the knot obtained from $K$ in
Figure~\ref{fig:alg-slice-with-nonslice-bd} by replacing the $\pm 3$
full twists on the leftmost and rightmost bands with $\pm q$ full
twists satisfies the conclusion of Theorem~\ref{thm:main3}.  (More
examples are given in Section~\ref{subsec:concrete}.)

We remark that the method in \cite{Cochran-Harvey-Leidy:2007-1} does
not give an explicit single knot with this property; using their
result, one can construct a family of knots such that all but possibly
one in the family should have the desired property, but their method
does not show which ones have the property.

It is known that the subgroup of algebraically slice knots in the
knot concordance group has a very rich structure. In
\cite{Cochran-Orr-Teichner:2002-1}, Cochran, Orr, and Teichner
constructed a filtration of the knot concordance group $\C$,
\[
0\subset \cdots\subset \F_{(n.5)}\subset \F_{(n)} \subset
\cdots\subset \F_{(1.5)}\subset \F_{(1)} \subset \F_{(0.5)}\subset
\F_{(0)}\subset \C,
\]
where $\F_{(h)}$ is the subgroup of \emph{$(h)$--solvable} knots. The
subgroup of algebraically slice knots is exactly $\F_{(0.5)}$, the
subgroup of $(0.5)$--solvable knots. Regarding this filtration, our
covering link method can also be used to produce examples which
satisfy the following refined statement:

\begin{theorem}[\cite{Cochran-Harvey-Leidy:2007-1}]
  \label{thm:refinement}
  For any integer $h\ge 1$, there are $(h)$--solvable knots $K$ such
  that for any $n$, $BD_n(K)$ is not slice.
\end{theorem}

Our examples and proofs are different from those given
in~\cite{Cochran-Harvey-Leidy:2007-1}.  To prove
Theorems~\ref{thm:main3} and~\ref{thm:refinement}, appealing to
Proposition~\ref{prop:reverse} stated above, it suffices to find an
algebraically slice or $(h)$--solvable knot $K$ for which $2K\# 2K^r$
is not rationally slice.  For this purpose we use von Neumann
$\rho$--invariants, which were used in \cite{Cha:2003-1} to give an
obstruction for algebraically slice knots to being rationally slice
(and to being linearly independent in the rational knot concordance
group).  For the highly solvable case of Theorem~\ref{thm:refinement},
we show that the examples in \cite{Cochran-Kim:2004-1} satisfy our
rational non-slice condition of $2K\# 2K^r$.  For this purpose, in
Section~\ref{sec:rational-concordance} we generalize some results on
integral knot concordance in \cite{Cochran-Kim:2004-1} to the rational
case.  Some arguments are essentially the same as the ones in
\cite{Cochran-Kim:2004-1} but some results in
Section~\ref{sec:rational-concordance} are not immediate consequences
of \cite{Cochran-Kim:2004-1}.  (Probably Theorem~\ref{thm:nontrivial}
and Proposition~\ref{prop:algebraic} are of independent interest.)

In fact using this approach we show a further generalization of
Theorem~\ref{thm:refinement}: there are highly solvable knots $K$
whose iterated Bing doubles are not only non-slice but also
non-solvable.  (For a precise statement, refer to
Theorem~\ref{thm:further-refinement}.)  For this purpose we use a
previous result of the first author called Covering Solution Theorem
\cite[Theorem~3.5]{Cha:2007-8}, which estimates solvability of
covering links.

As well as the above results that hold in both topological and smooth
cases, our covering link calculus method also gives results peculiar
to the smooth case: using Proposition~\ref{prop:reverse}, we show that
if $BD_n(K)$ is smoothly slice for some $n\ge 0$, then the
Heegaard--Floer homology theoretic concordance invariants of
Ozsv\'ath--Szab\'o \cite{Ozsvath-Szabo:2003-1} and Manolescu--Owens
\cite{Manolescu-Owens:2007-1} of $K$ vanish (see
Theorem~\ref{thm:heegaard-floer}).  This generalizes the special case
of $n=1$ proved in~\cite{Cha-Livingston-Ruberman:2006-1}.

The paper is organized as follows. In Section~\ref{sec:covering} we
define $p$--covering links and show their properties. We prove
Theorem~\ref{thm:main1} and Proposition~\ref{prop:reverse} in
Section~\ref{sec:zp-sliceness} and Theorem~\ref{thm:main2} in
Section~\ref{sec:algebraic}.  Our results on the Heegaard--Floer
invariants is proved in Section~\ref{sec:heegaard-floer}.
Theorem~\ref{thm:main3} and their refinements are proved in
Section~\ref{sec:rho-invariant}, and in
Section~\ref{sec:rational-concordance} we investigate rational knot
concordance and von Neumann $\rho$--invariants.

\subsection*{Acknowledgements}
The authors thank an anonymous referee for helpful comments.  This
work was supported by the Korea Science and Engineering
Foundation(KOSEF) grant funded by the Korea government(MOST) (No.\
R01--2007--000--11687--0).

\section{Covering links}\label{sec:covering}

Let $p$ be a prime and $\Sigma$ a $\zp$--homology 3--sphere. Note that a
manifold is a $\zp$--homology sphere if and only if it is a
$\Z_p$--homology sphere (and this is equivalent to that it is a
$\Z_{p^a}$-homology sphere for all/some~$a$).  Given a link $L$ in
$\Sigma$, we think of the following two operations producing new links
from~$L$:

\begin{enumerate}
\item[(C1)] Taking a sublink of $L$, a link in the same ambient space
  $\Sigma$ is obtained.

\item[(C2)] Choose a component $K$ of $L$ and a positive integer~$a$.
  From the homology long exact sequence for $(\Sigma,\Sigma-K)$ with
  $\Z_{p^a}$-coefficients and Alexander duality, we have
  \[
  H_1(\Sigma- K;\Z_{p^a})\cong H_2(\Sigma,\Sigma- K;\Z_{p^a}) \cong
  H^1(K;\Z_{p^a}) \cong \Z_{p^a}.
  \]
  Therefore there is a canonical map $\phi\colon H_1(\Sigma- K) \to
  \Z_{p^a}$ sending a meridian to a generator.  If the
  $(\Q/\Z)$--valued self-linking of $K$ in $\Sigma$ is trivial, then
  there is a ``preferred longitude'' of $K$ which is mapped to zero
  under the map $\phi$, due to~\cite{Cha:2003-1}. Therefore in this
  case the $p^a$--fold cyclic branched cover, say $\widetilde{\Sigma}$,
  of $\Sigma$ branched along $K$ is defined. By results of
  \cite{Casson-Gordon:1986-1} or more generally of
  \cite{Levine:1994-1}, $\widetilde{\Sigma}$ is a $\zp$--homology
  sphere and the preimage of $L$ can be viewed as a new link
  in~$\widetilde{\Sigma}$.
\end{enumerate}

\begin{definition}\label{definition:covering-link}
  A link $\widetilde{L}$ obtained from $L$ by applying (C1) and/or
  (C2) above repeatedly is called a \emph{$p$--covering link of $L$ of
    height $\le h$}, where $h$ is the number of (C2) applied.
\end{definition}

We remark that a different exponent $a$ can be used for each (C2)
applied.  As an abuse of terminology, we will often say that
$\widetilde L$ in Definition~\ref{definition:covering-link} is of
\emph{height~$h$}, although the precise definition of the height
should be the minimal number of (C2) applied.

It can be seen easily that if $L$ is a link in~$S^3$, the
$(\Q/\Z)$--linking number condition in (C2) above is automatically
satisfied.  Moreover, if a component $K$ of $L$ in a $\zp$--homology
sphere satisfies the condition as in (C2), then the condition also
holds for any component of the pre-image of $L$ in $\tilde\Sigma$
which projects to~$K$; for, due to \cite{Cha:2003-1}, $K$ satisfies
the $(\Q/\Z)$--linking number condition if and only if there is a
``generalized Seifert surface'' $F$, namely, an embedded oriented
surface $F$ in $\Sigma$ which is bounded by the union of $c>0$
parallel copies of $K$ taken along the zero-framing.  Considering a
component of the pre-image of $F$ in $\tilde\Sigma$, the claim easily
follows.  These observations enable us to iterate (C2) in
Definition~\ref{definition:covering-link} above in many cases.

Using the following well-known fact, we investigate the sliceness of a
link via its $p$--covering links:

\begin{theorem}\label{thm:covering slice}
  Let $p$ be a prime and $L$ a link in a $\zp$--homology
  sphere~$\Sigma$.  If $L$ is $\zp$--slice, then any $p$--covering link
  of $L$ is $\zp$--slice.
\end{theorem}
\begin{proof}
  A sublink $L'$ of $L$ is obviously a $\zp$--slice link. Suppose $L$
  bounds slice disks in a $\zp$--homology 4--ball $W$. Let
  $\widetilde{L}$ be the preimage of $L$ in $\widetilde{\Sigma}$,
  where $\widetilde{\Sigma}$ is a $\zp$--homology sphere obtained by
  taking a $p^a$--fold cyclic branched cover of $\Sigma$ branched along
  a component of $L$, say $K$.  By taking a $p^a$--fold cyclic branched
  cover of $W$ branched along the slice disk for $K$ in $W$, we obtain
  a 4--manifold $\widetilde{W}$ such that $\widetilde{\Sigma}
  = \partial \widetilde{W}$.  Due to~\cite{Casson-Gordon:1986-1},
  $\widetilde W$~is a $\zp$--homology ball, and the preimages of the
  slice disks for $L$ are slice disks in $\widetilde{W}$
  for~$\widetilde L$.
\end{proof}

The following construction of covering links will play a crucial role
for our purpose.

\begin{lemma}\label{lem:covering}
  Suppose $p$ is a prime and let $L_0$, $L_1$, and $L_2$ be the links
  in $S^3$ illustrated in Figures~\ref{fig:lemmalink0},
  \ref{fig:lemmalink1}, and \ref{fig:lemmalink2}, respectively. Then
  the following conclusions hold:
  \begin{enumerate}
  \item $L_1$ is a $p$--covering link of $L_0$ of height~1.
  \item $L_2$ is a $p$--covering link of $L_0$ of height~2.
  \end{enumerate}

  \begin{figure}[H]
    \begin{center}
      \includegraphics[scale=0.77]{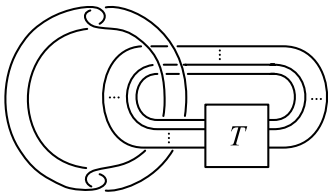}
    \end{center}
    \caption{Link $L_0$} \label{fig:lemmalink0}
  \end{figure}

  \begin{figure}[H]
    \begin{center}
      \includegraphics[scale=0.7]{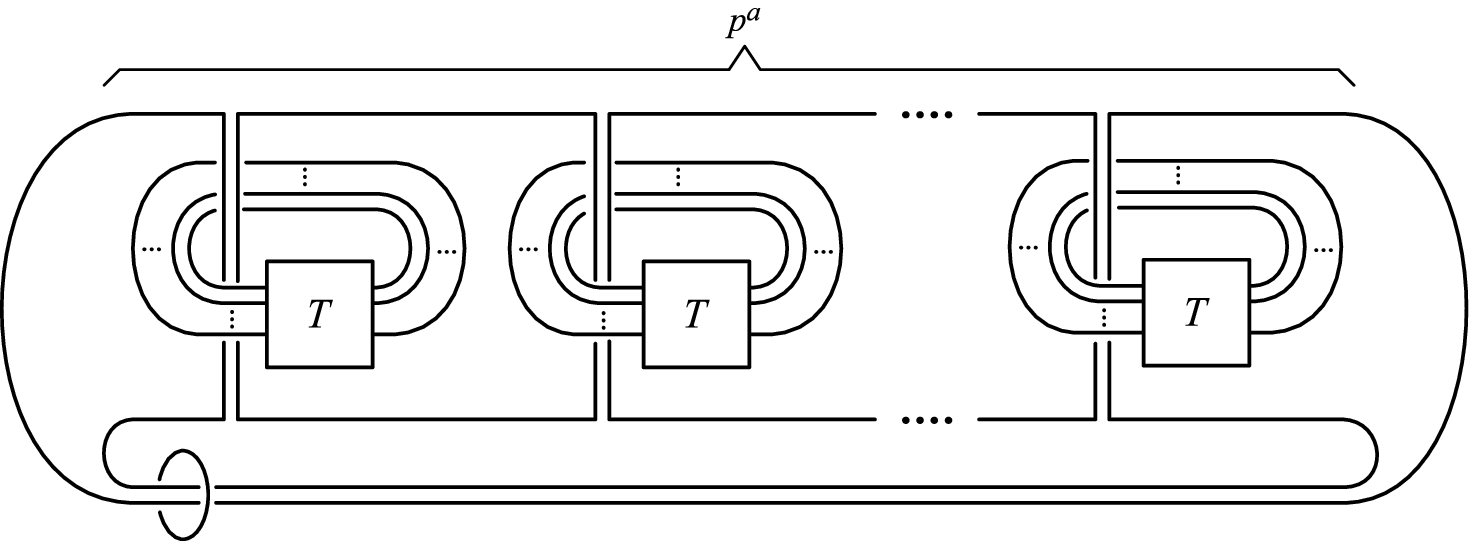}
    \end{center}
    \caption{Link $L_1$} \label{fig:lemmalink1}
  \end{figure}

  \begin{figure}[H]
    \begin{center}
      \includegraphics[scale=0.8]{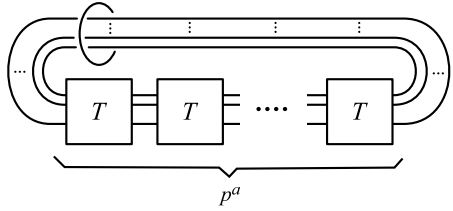}
    \end{center}
    \caption{Link $L_2$} \label{fig:lemmalink2}
  \end{figure}

\end{lemma}

\begin{figure}[H]
  \begin{center}
    \includegraphics[scale=0.77]{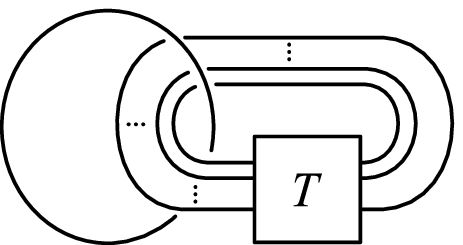}
  \end{center}
  \caption{Link $L_1'$} \label{fig:lemmalink1-1}
\end{figure}

\begin{proof}
  (1) The link $L_1$ is obtained by taking the $p^a$--fold cyclic
  branched cover of $S^3$ along the leftmost component of~$L_0$.

  (2) Forgetting appropriate components of $L_1$, we obtain the link
  $L_1'$ in Figure~\ref{fig:lemmalink1-1}.  Taking the $p^a$--fold
  cyclic branched cover of $S^3$ branched along the leftmost component
  of $L_1'$, we obtain~$L_2$.
\end{proof}

\section{Covering link construction relating iterated Bing doubles}
\label{sec:zp-sliceness}

\begin{figure}[H]
  \begin{center}
    \includegraphics[scale=0.8]{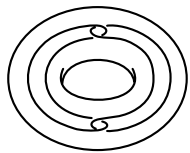}
  \end{center}
  \caption{The link $BD$ in a solid torus} \label{fig:bd}
\end{figure}

For clarity, we describe how the iterated Bing doubles are constructed
and fix notations.  In what follows a solid torus is always embedded
in $S^3$, so that its preferred longitude is defined.  Let $BD$ be the
2--component link contained in an unknotted solid torus illustrated in
Figure~\ref{fig:bd}.  For a link $L$, we define the Bing double
$BD_1(L)=BD(L)$ to be the link $L$ obtained by replacing a tubular
neighborhood of each component with a solid torus containing $BD$ in
such a way that a preferred longitude and a meridian of the solid
torus for $BD$ are matched up with those of the component of~$L$.  The
$n$th iterated Bing double $BD_n(L)$ is defined to be
$BD_n(L)=BD(BD_{n-1}(L))$.  For convenience, we denote $BD_0(L)=L$.

\begin{figure}[ht]
  \begin{center}
    \includegraphics[scale=0.8]{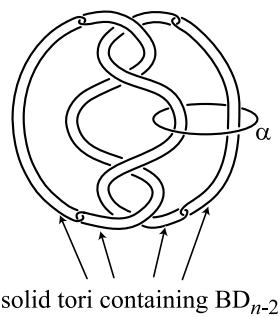}
  \end{center}
  \caption{$BD_n \cup \alpha$} \label{fig:alpha}
\end{figure}

For a knot $K$, we can construct $BD_n(K)$ using the process called
\emph{infection}.  A precise description is as follows.  Fix an
unknotted solid torus $V$ in $S^3$, and let $BD_0$ be the core of~$V$.
Let $BD_n=BD_n(BD_0)$ be the $2^n$--component link in~$V$.  Let
$\alpha$ be a meridional curve of~$V$.  See Figure~\ref{fig:alpha} for a
picture of $BD_n \cup\alpha$.
 As a simple closed curve in $S^3$, $\alpha$ is unknotted.
We take the union of the exterior of $\alpha\subset S^3$ and that of
the given knot $K\subset S^3$, glued along the boundaries such that a
longitude and a meridian for $\alpha$ are identified with a meridian
and a longitude for~$K$, respectively.  Then the resulting manifold is
homeomorphic to $S^3$, and $BD_n(K)$ is the image of $BD_n$ in this
new ambient manifold.

\begin{proposition}\label{prop:covering}
  Let $K$ be a knot in $S^3$. For any prime $p$ and any $n\ge 3$,
  $BD_{n-1}(K)$ is a $p$--covering link of $BD_n(K)$ of height 2.
\end{proposition}

From Proposition~\ref{prop:covering} and Theorem~\ref{thm:covering
  slice}, Theorem~\ref{thm:main1} follows immediately.

\begin{proof}[Proof of Proposition~\ref{prop:covering}]
  We regard $BD_n \cup \alpha$ as a link in $S^3$, and will show that
  $BD_{n-1}\cup \alpha$ is a $p$--covering link of $BD_n\cup \alpha$ by
  constructing a sequence of (C1) and (C2) operations.  In addition,
  we will observe that these operations behave in such a way that by
  performing infection along (pre-images of) $\alpha$, it follows that
  $BD_{n-1}(K)$ is a $p$--covering link of $BD_n(K)$ for any knot~$K$.

  We define $V_k$, $1\le k \le n$, to be a link in an unknotted solid
  torus as follows: first, $V_n$ is the core denoted by $\alpha$, as
  in the left in Figure~\ref{fig:v_k}.  For $k\le n$, we inductively
  define $V_{k-1}$ as in the right in Figure~\ref{fig:v_k}.  Note that
  each $V_k$ has a component denoted by~$\alpha$.

  \begin{figure}[H]
    \begin{center}
      \includegraphics[scale=0.8]{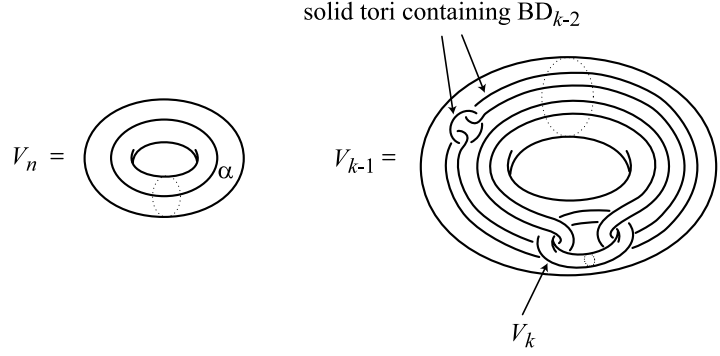}
    \end{center}
    \caption{$V_k$ for  $1\le k\le n$} \label{fig:v_k}
  \end{figure}

  $BD_n\cup \alpha$ can be illustrated as in the left in
  Figure~\ref{fig:v_n-1}.  (For convenience, the solid torus labeled
  as $V_k$ represents our link $V_k$ contained in the solid torus.)
  It can be seen that this is isotopic to the right diagram in
  Figure~\ref{fig:v_n-1}.  Note that we may denote this diagram by
  $BD_{n-1}\cup V_{n-1}$, by comparing it with Figure~\ref{fig:alpha}.

  \begin{figure}[H]
    \begin{center}
      \includegraphics[scale=0.8]{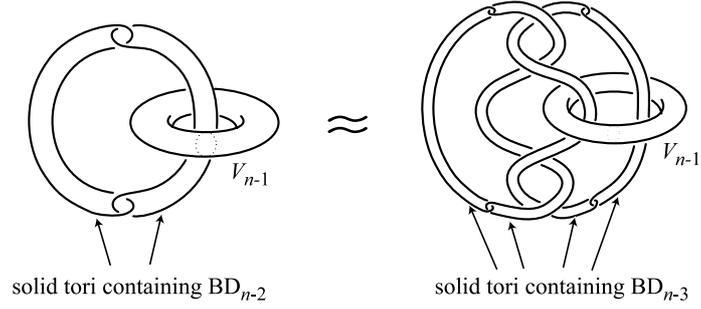}
    \end{center}
    \caption{$BD_{n} \cup \alpha$ isotoped to $BD_{n-1}\cup
      V_{n-1}$} \label{fig:v_n-1}
  \end{figure}

  Repeatedly applying this process, we have
  \[
  BD_n \cup \alpha \approx BD_n \cup V_n \approx BD_{n-1} \cup V_{n-1}
  \approx \cdots BD_1 \cup V_1
  \]
  where $BD_1 \cup V_1$ is illustrated in Figure~\ref{fig:v_1}.

  \begin{figure}[H]
    \begin{center}
      \includegraphics[scale=0.8]{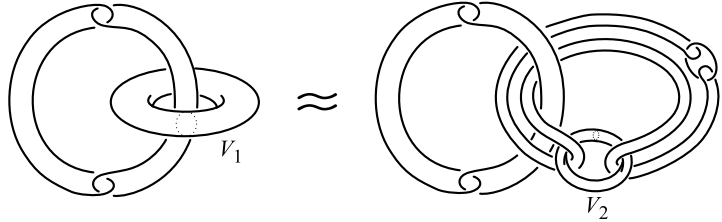}
    \end{center}
    \caption{$BD_{1}\cup V_{1}$} \label{fig:v_1}
  \end{figure}

  By Lemma~\ref{lem:covering}~(2), the link in
  Figure~\ref{fig:covering-of-bd_n} is a $p$--covering link of the link
  in the right in Figure~\ref{fig:v_1} (of height 2), hence of
  $BD_n\cup \alpha$.

  \begin{figure}[ht]
    \begin{center}
      \includegraphics[scale=0.8]{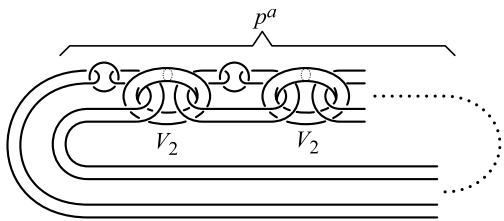}
    \end{center}
    \caption{A $p$--covering link of $BD_n\cup \alpha$}
    \label{fig:covering-of-bd_n}
  \end{figure}

  Forgetting some components of the link in
  Figure~\ref{fig:covering-of-bd_n}, we obtain the link in
  Figure~\ref{fig:v_2}.

  \begin{figure}[ht]
    \begin{center}
      \includegraphics[scale=0.8]{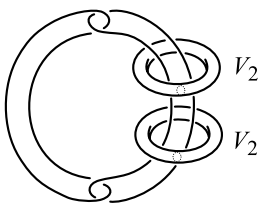}
    \end{center}
    \caption{Another $p$--covering link of $BD_n\cup
      \alpha$} \label{fig:v_2}
  \end{figure}

  Furthermore, since $n\ge 3$, $V_2 \ne V_n = \alpha$.  Therefore we
  can forget all components of the link in (the solid torus for) the
  second copy of $V_2$ in Figure~\ref{fig:v_2}, in order to obtain
  $BD_1 \cup V_2$.  In order to be precise, we need to be more careful
  with the component labeled as $\alpha$ in the second copy of $V_2$,
  since it is used as an infection curve.  Nonetheless, forgetting all
  components in the second $V_2$ but the concerned $\alpha$, one
  completely splits the $\alpha$ from the other remaining components,
  so that infection along $\alpha$ changes nothing.  We also note that
  one could not eliminate the second copy of $V_2$ in
  Figure~\ref{fig:v_2} if $V_2$ were~$\alpha$.

  Now we have that $BD_1 \cup V_2$ as a $p$--covering link of $BD_n\cup
  \alpha$.  Performing isotopies which were described above, we obtain
  \[
  BD_1\cup V_2 \approx BD_2\cup V_3 \approx \cdots \approx
  BD_{n-1}\cup V_n = BD_{n-1}\cup \alpha.
  \]
  It follows that $BD_{n-1}\cup \alpha$ is a $p$--covering link of
  $BD_n\cup \alpha$.
\end{proof}

For $n=2$, the proof of Proposition~\ref{prop:covering} shows the
following proposition:

\begin{proposition}\label{prop:reverse2}
For any prime $p$, $BD(K\#K^r)$ is a $p$--covering link of $BD_2(K)$
of height~2.
\end{proposition}

\begin{proof}
  As in the proof of Proposition~\ref{prop:covering}, the link in
  Figure~\ref{fig:v_2} is a $p$--covering link of $BD_2\cup \alpha$ of
  height 2.  Since $n=2$, one sees that $V_2 = \alpha$.  By carefully
  following the transform from Figure~\ref{fig:covering-of-bd_n} to
  Figure~\ref{fig:v_2}, one can see that the two copies of $\alpha
  (=V_2)$ in Figure~\ref{fig:v_2} are with opposite string
  orientations.  Performing infection by $K$ along the two copies of
  $\alpha$, the proposition follows.
\end{proof}

\begin{figure}[H]
  \begin{center}
    \includegraphics[scale=0.8]{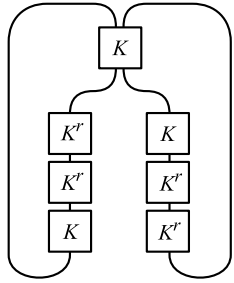}
  \end{center}
  \caption{A $p$--covering link $\widetilde{L}$ of $BD_2(K)$.}
  \label{fig:connected-sum}
\end{figure}

By arguments in \cite{Cha-Livingston-Ruberman:2006-1} or by applying
Lemma~\ref{lem:covering}~(1), it can be seen easily that $K\#K^r$ is a
$p$--covering link of $BD(K)$.  Consequently, by
Proposition~\ref{prop:reverse2}, the knot $2K\# 2K^r$ is a
$p$--covering link of $BD_2(K)$.  The following statement is a
generalization of this observation, which will be useful in
investigating algebraic invariants of iterated Bing doubles in
Section~\ref{sec:algebraic}:

\begin{proposition}\label{prop:connected sum}
Let $K$ be a knot in $S^3$. For every prime $p$, the link
$\widetilde{L}$ in Figure~\ref{fig:connected-sum} is a $p$--covering
link of $BD_2(K)$ of height 4.
\end{proposition}

\begin{proof}

  As in the proof of Theorem~\ref{thm:main1}, we start with the link
  $BD_2 \cup \alpha$ which is illustrated in
  Figure~\ref{fig:connected-sum-1}.

  \begin{figure}[H]
    \begin{center}
      \includegraphics[scale=0.8]{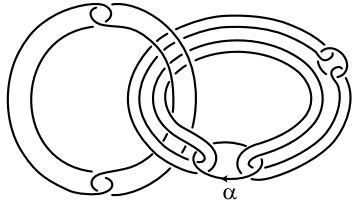}
    \end{center}
    \caption{Link $BD_2 \cup \alpha$} \label{fig:connected-sum-1}
  \end{figure}

  By Lemma~\ref{lem:covering}~(2), the link in the left in
  Figure~\ref{fig:connected-sum-2} is a $p$--covering link of $BD_2\cup
  \alpha$ of height~2.  Forgetting some components, we obtain the link
  in the right in Figure~\ref{fig:connected-sum-2}.

  \begin{figure}[H]
    \begin{center}
      \includegraphics[scale=0.8]{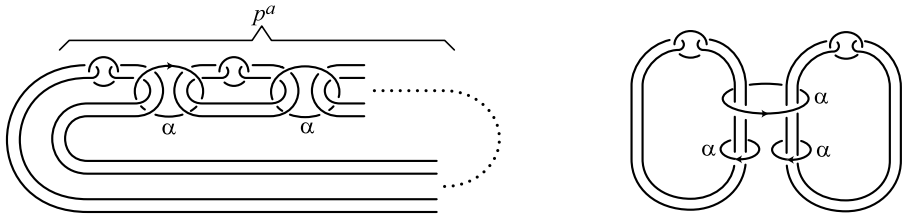}
    \end{center}
    \caption{$p$--covering links of $BD_2\cup \alpha$ of height 2}
    \label{fig:connected-sum-2}
  \end{figure}

  Applying Lemma~\ref{lem:covering}~(1), it follows that the link in
  the left in Figure~\ref{fig:connected-sum-3} is a $p$--covering link
  of $BD_2 \cup \alpha$ of height 3. Forgetting some components, we
  obtain the link in the right in Figure~\ref{fig:connected-sum-3}.

  \begin{figure}[H]
    \begin{center}
      \includegraphics[scale=0.8]{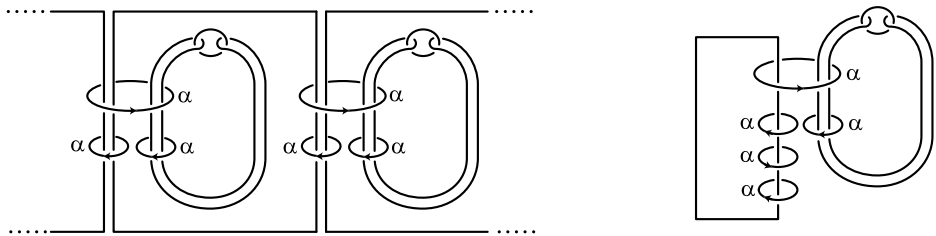}
    \end{center}
    \caption{$p$--covering links of $BD_2\cup \alpha$ of height 3}
    \label{fig:connected-sum-3}
  \end{figure}

  Again applying Lemma~\ref{lem:covering}~(1), it follows that the link
  in the left in Figure~\ref{fig:connected-sum-4} is a $p$--covering
  link of $BD_2 \cup \alpha$ of height 4. Forgetting some components,
  we obtain the link in the right in Figure~\ref{fig:connected-sum-4}
  as a $p$--covering link of $BD_2\cup \alpha$ of height 4.

\begin{figure}[H]
\begin{center}
\includegraphics[scale=0.8]{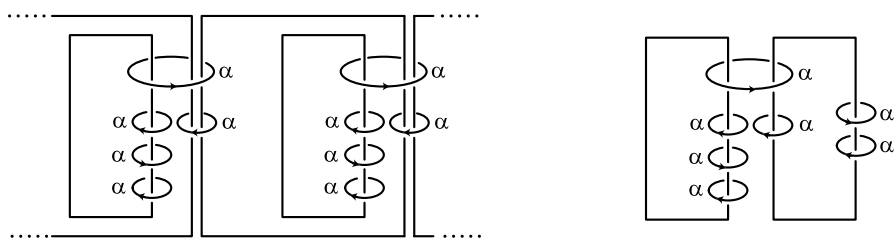}
\end{center}
\caption{$p$--covering links of $BD_2\cup \alpha$ of height 4}
\label{fig:connected-sum-4}
\end{figure}

Finally performing infection by $K$ along $\alpha$, it follows that
the link $\widetilde{L}$ in Figure~\ref{fig:connected-sum} is a
$p$--covering link of $BD_2(K)$ of height 4.
\end{proof}

\begin{proof}[Proof of Proposition~\ref{prop:reverse}]
  For $n=1$, the conclusion is known by arguments in
  \cite{Cha-Livingston-Ruberman:2006-1}.  (Or alternatively, apply
  Lemma~\ref{lem:covering}~(1) and Theorem~\ref{thm:covering slice}.)
  Suppose $n\ge 2$ and $BD_n(K)$ is $\zp$--slice.  By
  Theorem~\ref{thm:main1}, we may assume that $n=2$.  Note that
  $2K\#2K^r$ is a $p$--covering link of $BD_2(K)$ by forgetting one
  component of the link $\tilde L$ in Proposition~\ref{prop:connected
    sum}.  (Or alternatively, apply the paragraph above
  Proposition~\ref{prop:connected sum}.)  Therefore $2K\#2K^r$ is
  $\zp$--slice by Theorem~\ref{thm:covering slice}.
\end{proof}

\section{Algebraic invariants and iterated Bing doubles}
\label{sec:algebraic}

In this section we apply our geometric method to investigate algebraic
invariants of knots with slice iterated Bing doubles.  Recall that in
Section~\ref{sec:zp-sliceness} we showed that $2K \# 2K^r$ is
$\zp$--slice if $BD_n(K)$ is slice for some~$n$
(Proposition~\ref{prop:reverse}).  It can be seen that this conclusion
is strong enough to detect interesting examples of $K$ with non-slice
$BD_n(K)$ when the Levine--Tristram signature of $K$ is nontrivial, and
furthermore when a certain von Neumann $\rho$--invariant of $K$ is
nontrivial.  However, it gives no conclusion when $K$ is 2--torsion in
the (integral or $\zp$) knot concordance group, in particular when $K$
is amphichiral.  The first successful result on the non-sliceness of
$BD_n(K)$ for amphichiral $K$ was obtained in \cite{Cha:2007-5} using
invariants from iterated $p$--covers.  Our
Proposition~\ref{prop:connected sum} enables us to extract further
information when $K$ is amphichiral, via algebraic invariants of $K$,
as shown in the proof of Theorem~\ref{thm:main2} below.

\begin{proof}[Proof of Theorem~\ref{thm:main2}]
  Suppose $BD_n(K)$ is slice.  Let $A$ be a Seifert matrix of $K$ and
  $[A]$ be the element in the Levine's algebraic concordance group
  \cite{Levine:1969-2} represented by~$A$.  Our goal is to show that
  $[A]=0$.  For this purpose, we need the following facts on
  $\zp$--concordance: in \cite{Cha:2003-1} the algebraic
  $\zp$--concordance group and a canonical homomorphism from the
  algebraic concordance group to the algebraic $\zp$--concordance group
  are defined.  If a knot is $\zp$--slice, then its Seifert matrix
  represents a trivial element in the algebraic $\zp$--concordance
  group.  For $p=2$, it is known that the homomorphism of the
  algebraic concordance group to the algebraic $\z2$--concordance group
  is injective.  (For a detailed discussion on the necessary facts on
  $\zp$--concordance, see~\cite{Cha:2003-1}.)

  The map sending a knot $J$ to its $(c,1)$--cable $i_c(J)$ induces an
  endomorphism on the algebraic (integral and $\zp$) concordance
  group, and we denote the image of $[A]$ under this homomorphism by
  $i_c[A]$, following~\cite{Cha:2003-1}.  Consider the link $\tilde L$
  in Figure~\ref{fig:connected-sum}, which is $\zp$--slice by
  Propositions~\ref{prop:covering} and \ref{prop:connected sum} and
  Theorem~\ref{thm:covering slice}.  Taking one and $c$ parallel
  copies of the left and right components of $L$, respectively, and
  then attaching appropriate bands joining distinct components, we
  obtain a knot which is $\z2$--slice and has a Seifert matrix
  identical to that of the following connected sum:
  \[
  J=i_c(K)\,\,\#\,\, 2i_c(K^r)\,\, \#\,\, i_{c-1}(K)\,\, \#\,\, 2K\,\,
  \#\,\, K^r.
  \]
  Since $K$ and $K^r$ give the same element in the algebraic
  concordance group,
  \[
  3i_c[A] + i_{c-1}[A] + 3[A] = 0
  \]
  in the algebraic $\z2$ concordance group and thus in the algebraic
  concordance group.

  For $c=1$, we have $6[A]=0$. Since $4[A]=0$ whenever $[A]$ is
  torsion~\cite{Levine:1969-2}, it follows that $2[A]=0$. Therefore we
  have
  \[
  i_c[A] + i_{c-1}[A] + [A] = 0.
  \]
  By the arguments of \cite[Proof of Theorem
  1]{Cha-Livingston-Ruberman:2006-1}, it follows that $[A]=0$ in the
  (integral) algebraic concordance group.
\end{proof}

\section{Heegaard--Floer homology theoretic concordance invariants and iterated Bing doubles}
\label{sec:heegaard-floer}

In this section we consider two concordance invariants obtained from
Heegaard--Floer homology theory, namely the Ozsv\'ath--Szab\'o
$\tau$--invariant~\cite{Ozsvath-Szabo:2003-1} and Manolescu--Owens
$\delta$--invariant~\cite{Manolescu-Owens:2007-1}.  Our result can be
stated in a general form as follows:

\begin{theorem}
  \label{thm:heegaard-floer}
  Suppose $\phi$ is a torsion-free-abelian-group-valued knot invariant
  with the following properties:
  \begin{enumerate}
  \item $\phi$ is an invariant of unoriented knots, i.e.,
    $\phi(K)=\phi(K^r)$.
  \item $\phi$ is additive under connected sum, i.e.,
    $\phi(K_1\#K_2)=\phi(K_1)+\phi(K_2)$.
  \item $\phi$ is invariant under (smooth or topological)
    $\zp$--concordance for some prime $p$, i.e., $\phi(K)=0$ if $K$ is
    (smoothly or topologically) $\zp$--slice.
  \end{enumerate}
  If $BD_n(K)$ is (smoothly or topologically) slice for some $n$, then
  $\phi(K)=0$.
\end{theorem}

\begin{proof}
It follows immediately from Proposition~\ref{prop:reverse}.
\end{proof}

As mentioned in \cite[Section~4]{Cha-Livingston-Ruberman:2006-1},
$\tau$ and $\delta$ satisfy the above (1), (2), and (3) (for any $p$
and for $p=2$, respectively) in the smooth case.  Therefore, if
$BD_n(K)$ is smoothly slice for some $n$, then $\tau(K)=0$ and
$\delta(K)=0$.

\section{Von Neumann $\rho$--invariants and iterated Bing doubles}
\label{sec:rho-invariant}

In this section we construct algebraically slice knots with non-slice
iterated Bing doubles.  By Proposition~\ref{prop:reverse} the knot
$2K\#2K^r$ is $\zp$--slice for any prime~$p$ if $BD_n(K)$ is
$\zp$--slice for some~$n$.  Therefore for our purpose we will construct
algebraically slice knots $K$ such that $2K\#2K^r$ is not rationally
slice.

\subsection{Explicit examples}\label{subsec:concrete}

In \cite[Section 5]{Cha:2003-1}, it was shown that there exist
concrete and explicit examples of algebraically slice knots $K_i$,
$i\ge 1$, which are linearly independent in the rational knot
concordance group.  In particular, it was shown that for each $i$, the
knot $2K_i\# 2K_i (=4K_i)$ is not rationally slice.  Using the same
argument we will show that $2K_i\# 2K_i^r$ is not rationally slice, in
order to obtain the following theorem:

\begin{theorem}
  For the algebraically slice knots $K_i$ in \cite[Section
  5]{Cha:2003-1}, $BD_n(K_i)$ is not slice for any $n\ge 1$.
\end{theorem}

\begin{proof}
  First we describe the construction of $K=K_i$.  We choose a ``seed
  knot'' $K_0$ which is slice and has the rational Alexander module
  $\Q[t^{\pm 1}]/\langle p(t)^2 \rangle$ where $p(t)$ is a Laurent
  polynomial such that $p(t^{-1})$ equals $p(t)$ up to multiplication
  by a unit in $\Q[t^{\pm 1}]$, $p(1) = \pm 1$, and $p(t^c)$ is
  irreducible for any integer $c>0$.  The existence of such $p(t)$ and
  $K_0$ was shown in \cite[Section 5]{Cha:2003-1}. We choose a simple
  closed curve $\eta$ in $S^3- K_0$ which is unknotted in
  $S^3$ and satisfies the following: (1) $\ell k(\eta, K_0) = 0$, (2)
  the homology class $[\eta]$ in the rational Alexander module for
  $K_0$ equals $1 + \langle p(t)^2\rangle \in \Q[t^{\pm 1}]/\langle
  p(t)^2\rangle $. In particular, $[\eta]$ generates the rational
  Alexander module for~$K_0$.  Let $J$ be a knot such that $\rho(J)
  \ne 0$ where $\rho(J)$ denotes the integral of the Levine--Tristram
  signature function of $J$ over the unit circle normalized to length
  one. For example, one can take $J$ to be the connected sum of copies
  of the trefoil knot.  Then our $K=K_i$ is the knot $K_0(\eta, J)$,
  which denotes $K_0$ infected by a knot $J$ along the curve~$\eta$.

  In \cite[Proof of Theorem 5.25]{Cha:2003-1} it was shown that if
  $2K\# 2K'(\eta', J')$ is rationally slice for some slice knot $K'$
  and a simple closed curve $\eta'$, then
\[
\rho(J) + \epsilon \cdot \rho(J) + \epsilon' \cdot \rho(J')  = 0,
\]
for some nonnegative integers  $\epsilon \mbox{ and } \epsilon'$.
Note that $K^r = K_i^r = K_0^r(\eta, J)$. Therefore from the above
equation, if $2K\# 2K^r$ were rationally slice, then
\[
\rho(J) + \epsilon'' \cdot \rho(J)  = 0,
\]
for some nonnegative integer $\epsilon''$. Since $\rho(J)\ne 0$ by our
choice of $J$, one can conclude that $2K\# 2K^r$ is not rationally
slice, and the theorem follows.
\end{proof}

Theorem~\ref{thm:main3} stated in the introduction is a special case
of our construction in the above proof.  In fact, using a trefoil knot
as $J$, a genus 2 slice knot with $p(t)=3t^2-7t+3$ as $K_0$, we can
obtain the knot illustrated in
Figure~\ref{fig:alg-slice-with-nonslice-bd}.  The required
irreducibility condition is satisfied by Lemma~5.20
of~\cite{Cha:2003-1}.

\subsection{Examples of higher solvability}
\label{subsec:solvability}

In \cite{Cochran-Orr-Teichner:2002-1}, Cochran--Orr--Teichner
introduced a filtration on the knot concordance group $\C$
\[
0\subset \cdots\subset \F_{(n.5)}\subset \F_{(n)} \subset
\cdots \subset \F_{(1)} \subset \F_{(0.5)}\subset
\F_{(0)}\subset \C,
\]
where $\F_{(h)}$ is the subgroup of \emph{$(h)$--solvable} knots for
each nonnegative half-integer~$h$. The subgroup $\F_{(0.5)}$ is
exactly the subgroup of algebraically slice knots \cite[Remark
1.3.2]{Cochran-Orr-Teichner:2002-1}.
It is known that the filtration is nonstable.  For example,
$\F_{(h)}/\F_{(h.5)}$ is nontrivial for any integer~$h\ge
0$~\cite{Cochran-Orr-Teichner:1999-1,Cochran-Teichner:2003-1}.

We will show that for each integer $h>0$, there are $(h)$--solvable
knots $K$ such that $BD_n(K)$ is non-slice for any~$n$.  As we
discussed above, by Proposition~\ref{prop:reverse}, it suffices to
find $(h)$--solvable knots $K$ such that $2K\#2K^r$ is not rationally
slice.  We will show that certain examples of knots $K$ considered by
Cochran and the second author in \cite{Cochran-Kim:2004-1} have the
desired property.  In~\cite{Cochran-Kim:2004-1}, these knots $K$ were
shown to have the property that $4K=2K\# 2K$ is not (integrally)
slice, and it can be easily seen that their argument also shows that
$2K\#2K^r$ is not (integrally) slice.  We will show that it
generalizes to that $2K\#2K^r$ is not rationally slice.

Our argument is best described in terms of rational (or $\Z_{(p)}$-)
analogues of $(h)$--solvability and a relative version called
``rational $(h)$--solvequivalence'', which are due to
Cochran--Orr--Teichner \cite{Cochran-Orr-Teichner:2002-1}, Cha
\cite{Cha:2003-1,Cha:2007-8}, and Cochran--Kim
\cite{Cochran-Kim:2004-1}.  Although we will not use the definitions
excessively, we give precise definitions below for the convenience of
the reader.  For a group $G$, the \emph{$n$th derived group} $G^{(n)}$
of $G$ is defined inductively as follows: $G^{(0)}= G$ and
$G^{(n+1)}=[G^{(n)},G^{(n)}]$.  Let $R$ be a subring of~$\Q$.  The
examples to keep in mind are $R=\Z, \zp$, and $\Q$ where $p$ is a
prime.

\begin{definition}\label{def:solution}
    Let $n$ be a nonnegative integer and $W$ a 4--manifold with
    boundary components $M_1,\ldots,M_s$ such that $H_1(M_i;R)\cong R$
    for each~$i$.  The 4--manifold $W$ is called an
    \emph{$R$--coefficient $(n)$--cylinder} if the following hold:
    \begin{enumerate}
    \item $H_1(M_i;R) \rightarrow H_1(W;R)$ is an isomorphism for
      each~$i$, and
    \item there exist elements
      \[
      u_1, \ldots u_m, v_1,\ldots, v_m \in
      H_2\left(W;R[\pi/\pi^{(n)}]\right),
      \]
      where $m=\frac12 \dim_\Q \operatorname{Coker}\{ H_2(M;\Q) \to
      H_2(W;\Q)\}$ and $\pi=\pi_1(W)$, such that the
      $R[\pi/\pi^{(n)}]$-valued intersection form $\lambda_W^{(n)}$ on
      $H_2\big(W;R[\pi/\pi^{(n)}]\big)$ satisfies
      $\lambda_W^{(n)}(u_i,u_j) = 0$ and $\lambda_W^{(n)}(u_i,v_j) =
      \delta_{ij}$ (the Kronecker symbol).
    \end{enumerate}

  In addition, if the following holds then $W$ is called an
  \emph{$R$--coefficient $(n.5)$--cylinder}:
  \begin{itemize}
  \item[(3)] There exist $\tilde u_1, \ldots, \tilde u_m \in
    H_2\big(W;R[\pi/\pi^{(n+1)}]\big)$ such that
    $\lambda_W^{(n+1)}(\tilde u_i,\tilde u_j)=0$ and $u_i$ is the
    image of $\tilde u_i$ for each $i$.
  \end{itemize}
  Here the submodules generated by $\{u_i\}$, $\{v_i\}$, and $\{\tilde
  u_i\}$ are called an \emph{$(n)$--Lagrangian}, an \emph{$(n)$--dual},
  and an \emph{$(n+1)$--Lagrangian}, respectively.

  As a special case, for a nonnegative half-integer $h$, an
  $R$--coefficient $(h)$--cylinder $W$ with connected boundary $M$ is
  called an \emph{$R$--coefficient $(h)$--solution for~$M$}.
\end{definition}

\begin{definition}\label{def:solvability}
  Let $h$ be a nonnegative half-integer.  Two 3--manifolds $M$ and $M'$
  are \emph{$R$--coefficient $(h)$--solvequivalent} if there exists an
  $R$--coefficient $(h)$--cylinder $W$ such that $\partial W = M \coprod
  -M'$.  A 3--manifold $M$ is \emph{$R$--coefficient $(h)$--solvable} if
  there is a $R$--coefficient $(h)$--solution for~$M$.

  Two links in an $R$--homology 3--sphere are \emph{$R$--coefficient
    $(h)$--solvequivalent} if the zero surgeries on the links are
  $R$--coefficient $(h)$--solvequivalent.  A link in an $R$--homology
  3--sphere is \emph{$R$--coefficient $(h)$--solvable} if the zero
  surgery on the link is $R$--coefficient $(h)$--solvable.
\end{definition}

In the above definitions, when $R=\Q$, we often use ``rationally'' in
place of ``$\Q$--coefficient''.  Note that if a link is
($R$--coefficient) slice then it is $R$--coefficient $(h)$--solvable for
any subring $R$ of $\Q$ and for any~$h$.

In this subsection we only need a couple of facts on solvability and
solvequivalence (Proposition~\ref{prop:solvequivalece} and
Lemma~\ref{lem:observation}).  First, the following is an
$R$--coefficient version of \cite[Proposition 2.7]{Cochran-Kim:2004-1}.
The proof is identical to the argument in \cite{Cochran-Kim:2004-1},
and therefore we omit details.

\begin{proposition}\label{prop:solvequivalece}
  For two knots $J$ and $K$, if $J-K$ is $R$--coefficient
  $(h)$--solvable, then $J$ is $R$--coefficient $(h)$--solv\-equivalent
  to~$K$.  In particular, if $J-K$ is rationally slice, then $J$ is
  rationally $(h)$--solv\-equivalent to $K$ for all~$h$.
\end{proposition}

Here $-K$ is the mirror image with reversed orientation (i.e., a
concordance inverse) and $J-K$ denotes the connected sum of $J$
and~$-K$.

Denote the zero surgery manifold of a knot $K$ by~$M(K)$.  In
\cite{Cochran-Kim:2004-1}, for any given integer $h>0$, they
constructed an infinite family of certain knots $K_i$ such that for
any $i>j$, $K_i-K_j$ is $(h)$--solvable and $\coprod^k M(K_i)$ is not
$(h.5)$--solvequivalent to $\coprod^k M(K_j)$ whenever $k>0$.  The only
property of the $K_i$ we need is the following rational analogue:

\begin{lemma}\label{lem:observation}
  For any $i>j$ and $k>0$, $\coprod^k M(K_i)$ and $\coprod^k M(K_j)$
  are not rationally $(h.5)$--solvequivalent.
\end{lemma}

The proof of Lemma~\ref{lem:observation} is postponed to
Section~\ref{sec:rational-concordance}. (A precise description of
the $K_i$ is also given in Section~\ref{sec:rational-concordance}.)

\begin{proof}[Proof of Theorem~\ref{thm:refinement}]
  We will show that for the knot $K=K_i-K_j$ (with $i>j$), $BD_n(K)$
  is not slice for any~$n$. Since $ 2K\#2K^r = 2(K_i\#K_i^r)-
  2(K_j\#K_j^r) $, by Propositions~\ref{prop:reverse}
  and~\ref{prop:solvequivalece}, it suffices to show that
  $2(K_i\#K_i^r)$ is not rationally $(h.5)$--solv\-equivalent to
  $2(K_j\#K_j^r)$.

  Suppose that there is a rational $(h.5)$--cylinder, say $U$, between
  $M(2(K_i\#K_i^r))$ and $M(2(K_j\#K_j^r))$.  Note that for any finite
  collection $\{J_\ell\}$ of knots, there is a ``standard'' cobordism
  between $M(\#^\ell J_\ell)$ and $\coprod^\ell M(J_\ell)$ (e.g., see
  \cite[p.~113]{Cochran-Orr-Teichner:2002-1}).  Attaching such
  standard cobordisms to $U$, we obtain a rational $(h.5)$--cylinder
  between $\coprod^2 \big(M(K_i) \coprod M(K_i^r)\big)$ and $\coprod^2
  \big(M(K_j) \coprod M(K_j^r)\big)$.  Since $M(J)=M(J^r)$ for any
  knot $J$, we have actually obtained a rational $(h.5)$--cylinder
  betweeen $\coprod^4 M(K_i)$ and $\coprod^4 M(K_j)$.  This
  contradicts Lemma~\ref{lem:observation}.
\end{proof}

\subsection{Further refinement}\label{subsec:further-refinement}

In this subsection, we investigate relationships between
$\zp$--coefficient solvability of a link $L$ and that of a $p$--covering
link of $L$. The first interesting result along this line is the
\emph{Covering Solution Theorem} obtained by the first author in
\cite[Theorem 3.5]{Cha:2007-8}. For a space $X$ and a group
homomorphism $\pi_1(X)\rightarrow \Gamma$, let $X_\Gamma$ denote the
induced $\Gamma$--cover of~$X$.

\begin{theorem}\cite[Covering Solution Theorem]{Cha:2007-8}
  \label{thm:covering-solution-theorem}
  Let $p$ be a prime and $h\ge 1$ be a half-integer. Let $M$ be a
  closed 3--manifold. Suppose $W$ is a $\zp$--coefficient $(h)$--solution
  for $M$, $\phi\colon \pi_1(M) \rightarrow \Gamma$ is a homomorphism onto
  an abelian $p$--group $\Gamma$, and both $H_1(M)$ and $H_1(M_\Gamma)$
  are $p$--torsion free. Then $\phi$ extends to $\pi_1(W)$, and
  $W_\Gamma$ is an $(h-1)$--solution for $M_\Gamma$.
\end{theorem}

It immediately follows that (C2) in
Definition~\ref{definition:covering-link} reduces solvability of a
link by at most one:

\begin{corollary}\label{cor:covering-solution-theorem}
  Let $p$ be a prime and $h$ a half-integer with $h\ge 1$. Suppose $L$
  is a $\zp$--coefficient $(h)$--solvable link in a $\zp$--homology
  3--sphere and $\widetilde{L}$ is a $p$--covering link of $L$ obtained
  by applying (C2) in Definition~\ref{definition:covering-link} once.
  Then $\widetilde{L}$ is $\zp$--coefficient $(h-1)$--solvable.
\end{corollary}

On the other hand, the following theorem and its corollary show that
(C1) in Definition~\ref{definition:covering-link} preserves
solvability of a link.

\begin{theorem}\label{thm:sublink-solution}
  Let $M$ be a closed 3--manifold and $h$ a nonnegative half-integer.
  Suppose $W$ is an $R$--coefficient $(h)$--solution for $M$. Suppose
  $\alpha$ is a simple closed curve in $M$ such that the homology
  class $[\alpha]\in H_1(M;R)$ is of infinite order. Moreover, suppose
  that for the meridian $\mu_\alpha$ for $\alpha$, the homology class
  $[\mu_\alpha]=0$ in $H_1(M- \alpha;R)$. If $M'$ is the
  manifold obtained by surgery on $M$ along (any framing of) $\alpha$,
  then $M'$ is $R$--coefficient $(h)$--solvable.
\end{theorem}

\begin{proof}
Suppose $h$ is an integer. Let $W'$ be the manifold obtained from
$W$ by attaching a 2--handle along $\alpha$ in $M$. Then $\partial
W'=M'$ and
\begin{align*}
  H_1(W';R) & \cong  H_1(W;R)/\langle \alpha\rangle \cong
  H_1(M;R)/\langle \alpha\rangle,\\
  H_1(M;R) & \cong  H_1(M- \alpha;R)/\langle\mu_\alpha\rangle
  \cong  H_1(M- \alpha;R),\\
  H_1(M';R) & \cong
  H_1(M-\alpha;R)/\langle\lambda_\alpha\rangle \cong
  H_1(M;R)/\langle\alpha\rangle ,
\end{align*}
where $\lambda_\alpha$ is the longitude for $\alpha$. Therefore
$H_1(M';R) \rightarrow H_1(W';R)$ is an isomorphism. By
Mayer--Vietoris, we have the exact sequence
\[
0\to H_2(W;R) \to H_2(W';R) \to H_1(S^1;R)
\xrightarrow{i_*} H_1(W;R).
\]
Since $[\alpha]$ generates $H_1(S^1;R)$ and it is of infinite order in
$H_1(M;R)\cong H_1(W;R)$, the map $i_*$ is injective. It follows that
$H_2(W;R)\cong H_2(W';R)$. Therefore the images of the
$(h)$--Lagrangian and its $(h)$--dual for $W$ are an $(h)$--Lagrangian
and its $(h)$--dual for~$W'$.  Hence $W'$ is an $R$--coefficient
$(h)$--solution for $M'$.  When $h$ is a non-integral half-integer, the
theorem is similarly proved.
\end{proof}

\begin{corollary}\label{cor:sublink-solution}
  Suppose $L$ is an $R$--coefficient $(h)$--solvable link in an
  $R$--homology 3--sphere such that each component of $L$ has vanishing
  $R/\Z$--valued self linking number and any two distinct components
  have vanishing $R$--valued linking number.  Then, any sublink of $L$
  is $R$--coefficient $(h)$--solvable.
\end{corollary}

\begin{proof}
  It suffices to prove the theorem for the sublink $L' = L- K$
  where $K$ is a component of $L$. Let $M_L$ and $M_{L'}$ be the zero
  surgeries on $L$ and $L'$, respectively. Let $\alpha$ be the
  meridian for $K$.  Then $M_{L'}$ is homeomorphic to the manifold
  obtained from $M_L$ by surgery along~$\alpha$.

  Let $\mu_\alpha$ be the meridian for~$\alpha$.  From the
  self-linking number condition, it follows that there is a properly
  embedded oriented surface $F$ in the exterior of $K$ such that
  $\partial F$ is $c$ parallel copies of a preferred longitude of~$K$,
  where $c$ is an integer such that $1/c \in R$, due to
  \cite{Cha-Ko:2000-1}, \cite[Theorem 2.6(2)]{Cha:2003-1}.  (In
  \cite{Cha-Ko:2000-1}, \cite{Cha:2003-1}, such a surface $F$ is
  called a \emph{generalized Seifert surface for $K$ with complexity
    $c$}.)  Since the mutual linking number is zero, we may assume
  that $F$ is disjoint from $L - K$.  It follows that $c\mu_\alpha$ is
  homologous to $c\cdot($preferred longitude for $K)$, which is
  null-homologous in $M_L- \alpha$.  Thus $[\mu_\alpha]=0 $ in
  $H_1(M_L- \alpha;R)$.  Since $H_1(M_L;R)$ is freely generated by
  meridians for $L$, $[\alpha]$ is of infinite order in
  $H_1(M_L;R)$. Applying Theorem~\ref{thm:sublink-solution}, the
  corollary follows.
\end{proof}

\begin{corollary}\label{cor:p-covering-link-solution}
  Let $p$ be a prime and $h$ a nonnegative half-integer. Let $r$ be a
  nonnegative integer such that $r\le h$. Suppose $L$ is a
  $\zp$--coefficient $(h)$--solvable link in a $\zp$--homology 3--sphere
  and the linking number conditions in
  Corollary~\ref{cor:sublink-solution} are satisfied (here $R=\zp$).
  Then any $p$--covering link of $L$ of height $r$ is $\zp$--coefficient
  $(h-r)$--solvable.
\end{corollary}

\begin{proof}
It easily follows from
Corollaries~\ref{cor:covering-solution-theorem} and
\ref{cor:sublink-solution}.
\end{proof}

Using Corollary~\ref{cor:p-covering-link-solution} we can prove the
following theorem which strengthens Theorem~\ref{thm:refinement}. We
remark that Cochran, Harvey, and Leidy first proved (a more refined
version of) the following theorem in
\cite{Cochran-Harvey-Leidy:2007-1} using a different method and
examples.

\begin{theorem}\label{thm:further-refinement}
For any positive integers $h$ and $r$, there exists an
$(h)$--solvable knot $K$ such that $BD_r(K)$ is not
$(h+2r-0.5)$--solvable.
\end{theorem}

\begin{proof}
  Let $K$ be the knot $K_i-K_j$ considered in the proof of
  Theorem~\ref{thm:refinement}. Then $K$ is $(h)$--solvable. Suppose
  that $BD_r(K)$ is $(h+2r-0.5)$--solvable. Then it is
  $\zp$--coefficient $(h+2r-0.5)$--solvable. By
  Proposition~\ref{prop:covering}, $BD_2(K)$ is a $p$--covering link of
  $BD_r(K)$ of height $2r-4$. Then by Proposition~\ref{prop:reverse2}
  $BD(K\#K^r)$ is a $p$--covering link of $BD_r(K)$ of height $2r-2$.
  Recall that Lemma~\ref{lem:covering}(1) can be used to show that
  $J\#J^r$ is a $p$--covering link of $BD(J)$ of height~1.  Therefore
  it follows that $2K\#2K^r$ is a $p$--covering link of $BD_r(K)$ of
  height $2r-1$. By Corollary~\ref{cor:p-covering-link-solution} it
  follows that $2K\#2K^r$ is $\zp$--coefficient $(h.5)$--solvable. Since
  $\Q$ is flat over $\zp$, it follows that $2K\#2K^r$ is rationally
  $(h.5)$--solvable. But then as was shown in the proof of
  Theorem~\ref{thm:refinement}, it leads us to a contradiction.
\end{proof}

\section{Rational concordance and von Neumann $\rho$--invariants}
\label{sec:rational-concordance}

The purpose of this section is twofold: we extend results on
\emph{integral} concordance and solvability obtained by using the von
Neumann $\rho$--invariants in
\cite{Cochran-Teichner:2003-1,Cochran-Kim:2004-1} to the
\emph{rational} case, and give a proof of Lemma~\ref{lem:observation}
which was needed in the previous section.  If the reader is more
interested in the latter, we would recommend to read the last
subsection first, assuming Theorem~\ref{thm:eta}.

Essentially we follow the strategy of~\cite{Cochran-Kim:2004-1},
focusing on what differs from the integral case.  Details will be
omitted when arguments are almost identical to those of the integral
case.

\subsection{Homology of rational cylinders with PTFA coefficients}

To investigate rational $(n)$--cylinders more systematically, we need
the following notion of multiplicity given in \cite[Definition
2.1]{Cochran-Kim:2004-1}.  (It is often called the ``complexity'';
e.g., see \cite{Cochran-Orr:1993-1}, \cite{Cha-Ko:2000-1},
\cite{Cha:2003-1}.)

\begin{definition}
\label{def:multiplicity}
Let $h$ be a nonnegative half-integer. A boundary component $M$ of a
rational $(h)$--cylinder $W$ with $H_1(W;\Q)\cong \Q$ is said to be
\emph{of multiplicity $m$} if a generator in
$H_1(M)/\text{torsion}\cong \Z$ is sent to $m \in
H_1(W)/\text{torsion}\cong \Z$.
\end{definition}

We consider homology modules of 3--manifolds and rational cylinders
with coefficients in a certain Laurent polynomial ring
$\bbk[t^{\pm1}]$ over a skew field~$\bbk$, following the idea
of~\cite{Cochran-Orr-Teichner:1999-1} and subsequent works.  Details
are as follows.  Let $\Gamma$ be a poly-torsion-free-abelian
(henceforth PTFA) group such that $\Gamma/\G1 \cong \Z$, where
$\Gamma_r^{(n)}$ denotes the $n$th rational derived group of~$\Gamma$.
($\Gamma_r^{(0)}=\Gamma$ and $\Gamma_r^{(n)}$ is inductively defined
to be the minimal normal subgroup of $\Gamma_r^{(n-1)}$ such that
$\Gamma_r^{(n-1)}/\Gamma_r^{(n)}$ is abelian and torsion free; for
more details, see \cite[Section 3]{Harvey:2002-1}.)  Let $t$ be the
generator of~$\Z$.  Then $\Gamma \cong \Gamma^{(1)}_r\rtimes \langle
t\rangle$. Let $\K_\Gamma$ be the (skew) quotient field of
$\Z\Gamma$. The subgroup $\G1$ is also PTFA, and hence $\Z\G1$ embeds
in its (skew) quotient field, say~$\bbk$.  Therefore $\Z\Gamma =
\Z[\G1\rtimes \langle t\rangle]$ embeds in
$\Q\Gamma(\Q\G1-\{0\})^{-1}$, which is a Laurent polynomial ring
$\bbk[t^{\pm 1}]$.  Note $\Z\Gamma \subset \bbk[t^{\pm 1}] \subset
\K_\Gamma$ and $\bbk[t^{\pm 1}]$ is a PID.

Suppose $K$ is a knot with zero surgery $M$ and $W$ is a rational
$(n)$--cylinder which has $M$ as a boundary component of
multiplicity~$c$.  Suppose $\psi$ is a homomorphism of $\pi_1(W)$ into
our $\Gamma$ described above, which induces an isomorphism
$\pi_1(W)/\pi_1(W)^{(1)}_r \rightarrow \Gamma/\Gamma_r^{(1)} = \langle
t\rangle$.  Then $H_*(W;\bbk[t^{\pm 1}])$ is defined.  On the other
hand, the composition
\[
\pi_1(M) \to \pi_1(W)
\xrightarrow{\psi} \Gamma=\Gamma_r^{(1)} \rtimes \langle t \rangle
\]
factors through $\G1\rtimes \langle s\rangle $, where $s=t^c$.  As we
did for $\G1\rtimes \langle t\rangle $, $\Z[\G1\rtimes \langle
s\rangle]$ embeds into $\bbk[s^{\pm1}]$, and the homology module
$H_*(M;\bbk[s^{\pm 1}])$ is defined.  Viewing $\bbk[s^{\pm 1}]$ as a
subring of $\bbk[t^{\pm 1}]$, there is a natural map
\[
j_* \colon H_*(M;\bbk[s^{\pm 1}]) \to H_*(M;\bbk[t^{\pm 1}]).
\]

The following is a rational cylinder analogue of \cite[Theorem
3.8]{Cochran-Kim:2004-1}.

\begin{theorem}\label{thm:nontrivial}
  Suppose $K$, $M$, $W$, and $\Gamma$ are as above, and $\Gamma$ is
  $(n-1)$--solvable.  Let $d$ denote the degree of the Alexander
  polynomial of $K$. Then for
  the inclusion $i\colon M\hookrightarrow W$ we have
  \[
  \rank_\bbk \Im\left\{i_* \colon H_1(M;\bbk[t^{\pm 1}]) \rightarrow
    H_1(W;\bbk[t^{\pm 1}])\right\} \ge \left\{
    \begin{array}{ll}
      |c|(d-2)/2\,\, & \mbox{if } n>1, \\
      |c|d/2        & \mbox{if } n=1.
    \end{array} \right.
  \]
\end{theorem}

\begin{proof}
  Let $P = \Ker (i_*)$ and $Q = \Im (i_*)$. Let $\A' =
  H_1(M;\bbk[t^{\pm 1}])$ and $\A = H_1(M;\bbk[s^{\pm 1}])$.  It is
  known that the Blanchfield linking form
  $$
  \A' \to
  \Hom_{\bbk[t^{\pm 1}]}(\A', \K_\Gamma/\bbk[t^{\pm 1}])
  $$
  is nonsingular, and with respect to the Blanchfield linking form, $P
  \subset P^\perp$ \cite[Theorem 2.13]{Cochran-Orr-Teichner:1999-1}
  \cite[Proposition 3.6]{Cochran-Kim:2004-1}. Using these, one
  can show that
  $$
  \rank_\bbk Q \ge \frac12 \rank_\bbk
  \A',
  $$
  as done in the proof of \cite[Theorem
  3.8]{Cochran-Kim:2004-1}.

  Since $\bbk[s^{\pm 1}]$ is a (noncommutative) PID, we have a
  $\bbk[s^{\pm 1}]$-module isomorphism
  $$
  \A' \cong \A\mathbin{\mathop{\otimes}\limits_{\bbk[s^{\pm 1}]}}
  \bbk[t^{\pm 1}] \cong \bigoplus\limits^{|c|} \A
  $$
  as in \cite[Theorem 5.16(1)]{Cha:2003-1}.
  Therefore $\rank_\bbk \A' = |c| \cdot \rank_\bbk \A$, and it
  suffices to show that
  $$
  \rank_\bbk \A \ge \left\{\begin{array}{ll}
      d-2 & \mbox{ if } n >1, \\
      d   & \mbox{ if } n=1.
    \end{array} \right.
  $$
  
  Suppose $n>1$ and let $X=S^3-K$.  Since $\pi_1(X) \to
  \langle s\rangle $ is surjective, we have $H_1(X;\bbk[s^{\pm 1}])
  \cong H_1(X_\infty; \bbk)$, where $X_\infty$ denotes the connected
  infinite cyclic cover of~$X$. Therefore $\rank_\bbk
  H_1(X;\bbk[s^{\pm 1}]) \ge d-1$ by \cite[Corollary
  4.7]{Cochran:2002-1}.  Since the longitude for $K$ in
  $H_1(X;\bbk[s^{\pm 1}])$ is annihilated by $s-1\in \bbk[s^{\pm 1}]$
  and generates a $\bbk[s^{\pm 1}]$-submodule which is isomorphic to
  $\bbk$, we have $\rank_\bbk \A \ge d-2$.

  If $n=1$, $\Gamma$ is abelian and torsion free, and hence $\Gamma
  \cong \Z$. Therefore $\bbk = \Q$. Since $\pi_1(M)$ surjects to
  $\langle s\rangle$, $H_1(M;\Q[s^{\pm 1}])$ is the rational Alexander
  module.  It follows that $\rank_\bbk \A = d$.
\end{proof}

\begin{corollary}\label{cor:nontrivial}
  Suppose that $K$, $M$, $W$, and $\Gamma$ are as in
  Theorem~\ref{thm:nontrivial}. Let
  \[
  j_*\colon H_1(M;\bbk[s^{\pm 1}]) \to H_1(M;\bbk[t^{\pm 1}])
  \]
  be the map induced by the inclusion $j\colon \bbk[s^{\pm 1}]
  \to \bbk[t^{\pm 1}]$.  Then we have
  \[
  \rank_\bbk \Im\left\{i_* j_*\colon H_1(M;\bbk[s^{\pm 1}])
    \rightarrow H_1(W;\bbk[t^{\pm 1}])\right\} \ge \left\{
    \begin{array}{ll}
      (d-2)/2\,\, & \mbox{if } n>1, \\
      d/2        & \mbox{if } n=1.
    \end{array} \right.
  \]
\end{corollary}

\begin{proof}
As in the proof of Theorem~\ref{thm:nontrivial},
$$
H_1(M;\bbk[t^{\pm 1}]) \cong
\bigoplus\limits^{|c|} H_1(M;\bbk[s^{\pm 1}])
$$
as $\bbk[s^{\pm 1}]$-modules. The images of the $|c|$ copies of
$H_1(M;\bbk[s^{\pm 1}])$ under $i_*$ have the same
$\bbk$--rank since multiplication by $t^m$ $(m\in \Z)$ in
$H_1(W;\bbk[t^{\pm 1}])$ is an automorphism of
$H_1(W;\bbk[t^{\pm 1}])$ permuting the images of those
copies of $H_1(M;\bbk[s^{\pm 1}])$. Now the conclusion
follows from Theorem~\ref{thm:nontrivial}.
\end{proof}

\subsection{Rational cylinders and algebraic solutions}

In \cite{Cochran-Teichner:2003-1}, the notion of an algebraic
$(n)$--solution was first introduced in order to investigate the
behavior of $\pi_1(M) \to \pi_1(W) \to \pi_1(W)/\pi_1(W)_r^{(n)}$ for
a (integral) solution $W$ of~$M$.  In \cite{Cochran-Kim:2004-1},
Cochran and the second author extended it to (integral) cylinders.
For the convenience of the reader, the definition of an algebraic
$(n)$--solution \cite{Cochran-Kim:2004-1} is given below: for a group
$G$, let $G_k=G/G^{(k)}_r$. Then $G_k$ is a $(k-1)$--solvable PTFA
group, hence $\Z G_k$ embeds in its (skew) quotient field denoted by
$\bbk(G_k)$.

\begin{definition}
  \label{def:algebraic} Let S be a group such that $H_1(S;\Q)\neq 0$.
  Let $F$ be a free group and $i \colon F\to S$ a
  homomorphism. A homomorphism $r\colon S \to G$ is called an {\em
    algebraic $(n)$--solution ($n\ge0$) for $i \colon F\to S$}
  if the following hold:
  \begin{enumerate}
  \item For each $0\le k\le n-1$, the
    image of the following composition,  after tensoring with
    $\bbk(G_k)$, is nontrivial:
    $$
    H_1(S;\Z G_k)\overset{r_*}{\longrightarrow} H_1(G;\Z G_k)\cong
    G^{(k)}_{r}/[G^{(k)}_{r},G^{(k)}_{r}]\to
    G^{(k)}_{r}/G^{(k+1)}_{r}.
    $$
  \item For each $0\le k\le n$, the map $ H_1(F;\Z
    G_k)\overset{i_*}{\longrightarrow} H_1(S;\Z G_k)$, after tensoring
    with $\bbk(G_k)$, is surjective.
  \end{enumerate}
\end{definition}

The following is a generalization of \cite[Proposition
6.3]{Cochran-Kim:2004-1} to the case of \emph{rational}
$(n)$--cylinders:

\begin{proposition}\label{prop:algebraic}
  Suppose $n>0$ is an integer, $K$ is a knot with zero surgery~$M$,
  and the Alexander polynomial of $K$ has degree $d>2$.  (If $n=1$,
  $d=2$ is also allowed.)  Suppose that $W$ is a \emph{rational}
  $(n)$--cylinder with $M$ as one of its boundary components (of any
  multiplicity).  Let $\Sigma$ be a capped-off Seifert surface for~$K$.
  Suppose $F \to \pi_1(M- \Sigma)$ is a homomorphism of a free
  group $F$ inducing an isomorphism on $H_1(-;\Q)$.  Let
  $S=\pi_1(M)^{(1)}$, $G = \pi_1(W)^{(1)}_r$, and $i$ be the
  composition $F\to \pi_1(M- \Sigma) \to S$. Then the
  map $j\colon S \to G$ induced by inclusion is an algebraic
  $(n)$--solution for $i\colon F \to S$.
\end{proposition}

\begin{proof}
  We follow the lines in the proof of \cite[Proposition
  6.3]{Cochran-Kim:2004-1}. Let $\bbk = \bbk(G_k)$ be the (skew)
  quotient field of~$\Z G_k$.  First, we will prove that
  Definition~\ref{def:algebraic} (1) holds.  The map
  $G^{(k)}_r/[G^{(k)}_r, G^{(k)}_r] \rightarrow G^{(k)}_r/G^{(k+1)}_r$
  becomes an isomorphism after tensoring with $\bbk$, since its kernel
  is $\Z$--torsion.  Since $\bbk$ is flat over $\Z G_k$, it suffices to
  show that $j_*\colon H_1(S;\bbk) \rightarrow H_1(G;\bbk)$ is
  nontrivial.

  Let $c$ denote the multiplicity of $M$ for $W$. Let $\Gamma =
  \pi_1(W)/\pi_1(W)^{(k+1)}_r$. Then as in the previous subsection,
  $\Gamma \cong \Gamma/\G1 \rtimes \langle t \rangle$, the composition
  $\pi_1(M) \rightarrow \pi_1(W) \rightarrow \Gamma$ factors through
  $\Gamma/\G1 \rtimes \langle s\rangle$ where $s=t^c$, and
  $H_*(M;\bbk[s^{\pm 1}])$ and $H_*(W;\bbk[t^{\pm 1}])$ are defined.
  Since $\pi_1(M)/\pi_1(M)^{(1)} = \langle s\rangle$ and $S =
  \pi_1(M)^{(1)}$, we have $H_1(S;\bbk) \cong H_1(M;\bbk[s^{\pm 1}])$.
  Similarly, $H_1(G;\bbk) \cong H_1(W;\bbk[t^{\pm 1}])$.  Therefore
  $j_*$ is identical to $H_1(M;\bbk[s^{\pm 1}])\rightarrow
  H_1(W;\bbk[t^{\pm 1}])$.  By Corollary~\ref{cor:nontrivial}, it is
  nontrivial.

  One can prove that Definition~\ref{def:algebraic} (2) holds using
  the argument of the proof in \cite[Proposition
  6.3]{Cochran-Kim:2004-1}; one only needs to replace $\bbk[t^{\pm
    1}]$ by $\bbk[s^{\pm 1}]$.
\end{proof}

We have the following theorem which generalizes \cite[Theorem
5.13]{Cochran-Kim:2004-1} to the case of \emph{rational}
$(n)$--cylinders.

\begin{theorem}\label{thm:eta}
  Suppose $n$, $K$, and $M$ are as in
  Proposition~\ref{prop:algebraic}.  For any given Seifert surface for
  $K$, there exists an oriented trivial link $\{\eta_1, \eta_2,
  \ldots, \eta_m\}$ in $S^3$ which is disjoint from the Seifert
  surface and satisfies the following:
  \begin{enumerate}
  \item $\eta_i\in \pi_1(M)^{(1)}$ for all~$i$.  Furthermore, the
    $\eta_i$ bound (smoothly embedded) symmetric capped gropes of
    height $n$, disjointly embedded in $S^3- K$ (except for the caps,
    which may intersect~$K$).
  \item For every \emph{rational} $(n)$--cylinder $W$ with $M$ as one
    of its boundary components (of any multiplicity), there is some
    $\eta_i$ such that $j_*(\eta_i)\notin \pi_1(W)^{(n+1)}_r$ where
    $j_*\colon \pi_1(M) \rightarrow \pi_1(W)$ is induced by the
    inclusion. The number of such $\eta_i$ is at least $(d-2)/2$ if
    $n>1$ or at least $d/2$ if $n=1$ where $d$ denotes the degree of
    the Alexander polynomial for $K$.
  \end{enumerate}
\end{theorem}

\begin{proof}
  One can proceed exactly in the same way as the proof of the integral
  version \cite[Theorem 5.13]{Cochran-Kim:2004-1}, except that one
  should use $\pi_1(W)^{(1)}_r$, $H_1(M;\bbk(G_{n-1})[s^{\pm 1}])$ and
  Corollary~\ref{cor:nontrivial} and Proposition~\ref{prop:algebraic}
  proved above, instead of $\pi_1(W)^{(1)}$,
  $H_1(M;\allowbreak\bbk(G_{n-1})[t^{\pm 1}])$, and the integral analogues used
  in \cite{Cochran-Kim:2004-1}.  (Here $t$, $s$ are as in the proof of
  Proposition~\ref{prop:algebraic}.)
\end{proof}

We remark that the $\eta_i$ in Theorem~\ref{thm:eta} are the same as
those used in \cite[Theorem 5.13]{Cochran-Kim:2004-1}.

\subsection{Rational knot concordance and the Cochran--Orr--Teichner filtration}
\label{subsec:rational-knot-concordance-COT-filtration}

For a given postive integer $n$, we consider a family of knots $K_i$
which was given in \cite[Theorem 5.1]{Cochran-Kim:2004-1}.  For the
convenience of the reader, we briefly describe how the $K_i$ are
constructed.  Choose $K$ and $\{\eta_1,\ldots,\eta_m\}$ satisfying
(the conclusion of) Theorem~\ref{thm:eta}.  Then we let $K_0 = K$ and
for $i\ge 1$, $K_i = K(\eta_1, \ldots, \eta_m, J^i_1, \ldots, J^i_m)$,
the knot obtained from $K$ by infection along the $\eta_\ell$, where
the infection knots $J^i_\ell$ are chosen so that (the integrals of)
the Levine--Tristram sigantures of the $J^i_\ell$ satisfy certain
inequalities described in \cite[p.~1429, proof of
Theorem~5.1]{Cochran-Kim:2004-1}.  For more details, refer
to~\cite{Cochran-Kim:2004-1}.

For the $K_i$, we prove Lemma~\ref{lem:observation} used in the
previous section: $\coprod^k M(K_i)$ and $\coprod^k M(K_j)$ are not
\emph{rationally} $(n.5)$--solv\-equivalent.

\begin{proof}[Proof of Lemma~\ref{lem:observation}]
  We follow the arguments of the proof of \cite[Theorem
  5.1(5)]{Cochran-Kim:2004-1}, which shows that $\coprod^k M(K_i)$ and
  $\coprod^k M(K_j)$ are not \emph{integrally}
  $(n.5)$--solv\-equivalent.  All the arguments of
  \cite{Cochran-Kim:2004-1} proving their integral statement work
  verbatim in our case except that Theorem~\ref{thm:eta} should be
  applied instead of \cite[Theorem 5.13]{Cochran-Kim:2004-1}, in order
  to guarantee that whenever $W$ is a \emph{rational} $(n)$--cylinder
  with $M(K)$ as one of its boundary components, $j_*(\eta_\ell)\notin
  \pi_1(W)^{(n+1)}_r$ for some~$\eta_\ell$.
\end{proof}

Consider the Cochran--Orr--Teichner filtration
\[
0\subset \cdots\subset \F^\Q_{(n.5)}\subset \F^\Q_{(n)} \subset
\cdots\subset \F^\Q_{(1)} \subset \F^\Q_{(0.5)}\subset
\F^\Q_{(0)}\subset \C^\Q
\]
of the rational knot concordance group $\C^\Q$ \cite{Cha:2003-1},
where $\F^\Q_{(h)}$ is the subgroup of rationally $(h)$--solvable
knots.

\begin{theorem}\label{thm:r-solvable}
  For the $K_i$, the following hold:
  \begin{enumerate}
  \item If $i\ne j$, $K_i$ is not \emph{rationally}
    $(n.5)$--solvequivalent to $K_j$. In particular, $K_i - K_j$ is not
    \emph{rationally} $(n.5)$--solvable.
  \item For each $i> j$, $K_i-K_j$ is of infinite order in
    $\F^\Q_{(n)}/\F^\Q_{(n.5)}$.
  \end{enumerate}
\end{theorem}

\noindent The corollary below, which was first proved in \cite[Theorem
4.3]{Cochran-Harvey-Leidy:2007-1}, easily follows from
Theorem~\ref{thm:r-solvable}.  We remark that it is further
generalized to an infinite rank result
in~\cite{Cochran-Harvey-Leidy:2007-2}.

\begin{corollary}\label{cor:r-solvable}
For each positive integer $n$, $\F^\Q_{(n)}/\F^\Q_{(n.5)}$ has
positive rank.
\end{corollary}

\begin{proof}[Proof of Theorem~\ref{thm:r-solvable}]
  The first part is a special case of Lemma~\ref{lem:observation}
  (when $k=1$).  For the second part, suppose that for a positive
  integer $k$, the connected sum $\#^k(K_i - K_j)$ is
  $(n.5)$--solvable. Then by Proposition~\ref{prop:solvequivalece},
  $M(\#^k K_i)$ and $M(\#^k K_j)$ are $(n.5)$--solv\-equivalent. Let
  $U$ denote an $(n.5)$--cylinder between $M(\#^k K_i)$ and $M(\#^k
  K_j)$. As we did in the proof of Theorem~\ref{thm:refinement}, by
  attaching standard cobordisms to $U$ we obtain an $(h.5)$--cylinder
  between $\coprod^k M(K_i)$ and $\coprod^k M(K_j)$. This contradicts
  Lemma~\ref{lem:observation}.
\end{proof}

\bibliographystyle{amsplainabbrv}
\renewcommand{\MR}[1]{}

\bibliography{research}

\end{document}